\documentclass[12pt]{amsart}

 \usepackage{amsfonts,graphics,amsmath,amsthm,amsfonts,amscd,amssymb,amsmath,latexsym,multicol}
\usepackage{epsfig,url}
\usepackage{flafter}

\makeatletter

\def\jobis#1{FF\fi
  \def\predicate{#1}%
  \edef\predicate{\expandafter\strip@prefix\meaning\predicate}%
  \edef\job{\jobname}%
  \ifx\job\predicate
}

\makeatother

\if\jobis{proposal}%
\else
\fi

 \usepackage[matrix, arrow]{xy}

\DeclareMathOperator{\Supp}{Supp}

 \newcommand{\Q}{\mathbb Q}
 \newcommand{\R}{\mathbb R}
 
 \newcommand{\bir}{\dashrightarrow}
  \newcommand{\ep}{\varepsilon}

 \numberwithin{equation}{subsection}
 \numberwithin{footnote}{subsection}

 \newtheorem{thm}[subsection]{Theorem}

{
    \newtheoremstyle{upright}%
        {8pt plus2pt minus4pt}%
        {8pt plus2pt minus4pt}%
        {\upshape}%
        {}%
        {\bfseries\scshape}%
        {}%
        {1em}%
        {}%

\theoremstyle{upright}

 \newtheorem{defn}[subsection]{Definition}

}

 \newcommand{\ke}[1]{$\acute{\mbox{e}}$}
 \newcommand{\ku}[1]{$\acute{\mbox{u}$}}
 \newcommand{\kl}[1]{$\acute{\mbox{l}}$}
 \newcommand{\kh}[1]{$\acute{\mbox{h}}$}
 \newcommand{\kr}[1]{$\acute{\mbox{r}}$}
 \newcommand{\kx}[1]{$\acute{\mbox{x}}$}
 \newcommand{\ki}[1]{${\^\i}$}

\title{Supplement to the paper\\ "On existence of log minimal models II"}
\author{Caucher Birkar}\thanks{2000 Mathematics Subject Classification: 14E30}
\date{\today}
\begin{document}
\maketitle

\begin{abstract}
We prove a stronger version of a termination theorem appeared in the paper 
"On existence of log minimal models II"  [\ref{B-II}]. We essentially just get rid of 
the redundant assumptions so the proof is almost the same as in 
[\ref{B-II}]. However, we give a detailed proof here for future reference.
\end{abstract}


\section{Introduction}

We work over a fixed algebraically closed field $k$ of characteristic zero. See section $2$ 
of [\ref{B-II}] for notation, terminology, and definition and basic properties of log minimal models.
The following theorem was proved in [\ref{B-II}, Theorem 1.5] under stronger assumptions.
However, the proof works with little change under weaker assumptions.

\begin{thm}\label{t-main2}
 Let $(X/Z,B+C)$ be a lc pair of dimension $d$ such that $K_X+B+C$ is nef$/Z$, 
 $B,C\ge 0$ and $C$ is $\R$-Cartier.
Assume that we are given an LMMP$/Z$ on $K_X+B$ with scaling of $C$ as in 
Definition \ref{d-scaling} with $\lambda_i$ the corresponding numbers, 
and $\lambda:=\lim\lambda_i$. 
Then, the LMMP terminates in the following cases:

$(i)$ $(X/Z,B)$ is $\Q$-factorial dlt, $B\ge H\ge 0$ for some ample$/Z$ 
$\R$-divisor $H$;

$(ii)$ $(X/Z,B)$ is $\Q$-factorial dlt, 
$C\ge H\ge 0$ for some ample$/Z$ $\R$-divisor $H$, and $\lambda>0$;

$(iii)$ $(X/Z,B+\lambda C)$  has a log minimal model, and $\lambda\neq \lambda_j$ for any $j$.\\
\end{thm}

We should remark that much of the difficulties in the proof of  Theorem \ref{t-main2} are caused 
by the presence of non-klt singularities. It is also worth to mention that the above theorem and 
the corresponding Theorem 1.5 of [\ref{B-II}] follow quite different goals. The arguments of 
[\ref{B-II}]  in essence do not rely on [\ref{BCHM}]. In contrast, the above theorem heavily relies on 
[\ref{BCHM}]. So, this paper is not intended to replace [\ref{B-II}].

A \emph{sequence of log flips$/Z$ starting with} $(X/Z,B)$ is a sequence $X_i\bir X_{i+1}/Z_i$ in which  
$X_i\to Z_i \leftarrow X_{i+1}$ is a $K_{X_i}+B_i$-flip$/Z$, $B_i$ is the birational transform 
of $B_1$ on $X_1$, and $(X_1/Z,B_1)=(X/Z,B)$. As usual, here  $X_i\to Z_i$ is an extremal flipping 
contraction.

\begin{defn}[LMMP with scaling]\label{d-scaling}
Let $(X_1/Z,B_1+C_1)$ be a lc pair such that $K_{X_1}+B_1+C_1$ is nef/$Z$, $B_1\ge 0$, and $C_1\ge 0$ is $\R$-Cartier. 
Suppose that either $K_{X_1}+B_1$ is nef/$Z$ or there is an extremal ray $R_1/Z$ such
that $(K_{X_1}+B_1)\cdot R_1<0$ and $(K_{X_1}+B_1+\lambda_1 C_1)\cdot R_1=0$ where
$$
\lambda_1:=\inf \{t\ge 0~|~K_{X_1}+B_1+tC_1~~\mbox{is nef/$Z$}\}
$$
 When $(X_1/Z,B_1)$ is $\Q$-factorial dlt, the last sentence follows from [\ref{B}, 3.1] 
 (the same is true in general for lc pairs by the results of Ambro [\ref{Ambro}] and Fujino [\ref{Fujino}, Theorem 1.1 (6)]; however, we do not need this stronger version). 
If $R_1$ defines a Mori fibre structure, we stop. Otherwise assume that $R_1$ gives a divisorial 
contraction or a log flip $X_1\bir X_2$. We can now consider $(X_2/Z,B_2+\lambda_1 C_2)$  where $B_2+\lambda_1 C_2$ is 
the birational transform 
of $B_1+\lambda_1 C_1$ and continue. That is, suppose that either $K_{X_2}+B_2$ is nef/$Z$ or 
there is an extremal ray $R_2/Z$ such
that $(K_{X_2}+B_2)\cdot R_2<0$ and $(K_{X_2}+B_2+\lambda_2 C_2)\cdot R_2=0$ where
$$
\lambda_2:=\inf \{t\ge 0~|~K_{X_2}+B_2+tC_2~~\mbox{is nef/$Z$}\}
$$
 By continuing this process, we obtain a sequence of numbers $\lambda_i$ and a 
special kind of LMMP$/Z$ which is called the \emph{LMMP$/Z$ on $K_{X_1}+B_1$ with scaling of $C_1$}; note that it is not unique. 
This kind of LMMP was first used by Shokurov [\ref{log-flips}].
When we refer to \emph{termination with scaling} we mean termination of such an LMMP. We usually put 
$\lambda=\lim \lambda_i$.

When we have a lc pair $(X/Z,B)$, we can always find an ample$/Z$ $\R$-Cartier divisor $C\ge 0$ such that 
$K_X+B+C$ is lc and nef$/Z$,  so we can run the LMMP$/Z$ with scaling assuming that all the 
necessary ingredients exist, e.g. extremal rays, log flips. 
\end{defn}


\vspace{0.1cm}
\section{Proof of the theorem}

\vspace{0.3cm}
{\textbf{Proof of Theorem \ref{t-main2}}.}
\emph{Step 1.} Proof of (i):  Since $H$ 
is ample$/Z$, we can perturb the coefficients of $B$ hence assume that 
$(X/Z,B)$ is klt. If $\lambda_i<1$ for some $i$, then  
$({X}/Z,B+\lambda_iC)$ is klt. If $\lambda_i=1$ for every $i$, 
then we can perturb $C$ hence in any case we could assume that 
$(X/Z,B+C)$ is klt and then use [\ref{BCHM}]. 

Proof of (ii): The LMMP is also an LMMP$/Z$ on $K_{X}+B+\frac{\lambda}{2} C$
with scaling of $(1-\frac{\lambda}{2})C$. We can replace $B$ with $B+\frac{\lambda}{2} C$, 
$C$ with $(1-\frac{\lambda}{2})C$, and  $\lambda_i$ with 
$\frac{\lambda_i-\frac{\lambda}{2}}{1-\frac{\lambda}{2}}$. After this change, 
 we can assume that $B\ge \frac{H}{2}$. Now use (i).

Proof of (iii): Note that if\\

$(*)$ $(X/Z,B)$ is $\Q$-factorial dlt, and $C\ge H\ge 0$ for some ample$/Z$ $\R$-divisor $H$,\\\\
then, for each $i$, there is a klt $K_X+\Delta\sim_\R K_X+B+\lambda_i C/Z$. 
We continue the proof \emph{without} assuming $(*)$ (but we will come back to $(*)$ 
in Step 7).\\

\emph{Step 2.}  
We can replace $B$ with $B+\lambda C$ hence 
assume that $\lambda=0$. Moreover, we may assume that the LMMP consists of 
only a sequence $X_i\bir X_{i+1}/Z_i$ of log flips starting with 
$(X_1/Z,B_1)=(X/Z,B)$. Pick $i$ so that $\lambda_i>\lambda_{i+1}$. Thus, $\Supp C_{i+1}$ does not 
contain any lc centre of $(X_{i+1}/Z,B_{i+1}+\lambda_{i+1}C_{i+1})$ because  $(X_{i+1}/Z,B_{i+1}+\lambda_{i}C_{i+1})$ is 
lc. Then, by replacing $(X/Z,B)$ with $(X_{i+1}/Z,B_{i+1})$ 
and $C$ with $\lambda_{i+1}C_{i+1}$ we may 
assume that no lc centre of $(X/Z,B+C)$ is inside $\Supp C$. Moreover, since 
there are finitely many lc centres of $(X/Z,B)$, perhaps after truncating the 
sequence, we can assume that no lc centre is contracted in the 
sequence.

 By assumptions there is a log minimal model $(Y/Z,B_Y)$ for $(X/Z,B)$.
Let $\phi \colon X\bir Y/Z$ be the corresponding birational map. Since $K_{X_i}+B_i+\lambda_i C_i$ 
is nef$/Z$, we may add an ample$/Z$ $\R$-divisor $G^i$ so that $K_{X_i}+B_i+\lambda_i C_i+G^i$ 
becomes ample$/Z$, in particular, it is movable$/Z$. We can choose the $G^i$ so that 
$\lim_{i\to \infty} {G^i}_1=0$ in $N^1(X_1/Z)$ where ${G^i}_1$ is the birational transform of $G^i$ 
on $X_1=X$. Therefore, 
$$
K_{X}+B\equiv \lim_{i\to \infty} (K_{X_1}+B_1+\lambda_i C_1+G^i_1)/Z
$$ 
which implies that $K_{X}+B$ is a limit of movable$/Z$ $\R$-divisors.  

Let $f\colon W\to X$ and $g\colon W\to Y$ 
be a common log resolution of $(X/Z,B+C)$ and $(Y/Z,B_Y+C_Y)$ where $C_Y$ is the birational transform of $C$. 
By applying the negativity lemma to $f$, we see that 
$$
E:=f^*(K_X+B)-g^*(K_Y+B_Y)=\sum_D a(D,Y,B_Y)D-a(D,X,B)D
$$
is effective (cf. [\ref{B}, Remark 2.6]) where $D$ runs over the prime divisors on $W$. 
Assume that $D$ is a component of $E$. If $D$ is not exceptional$/Y$, 
then it must be exceptional$/X$ otherwise $a(D,X,B)=a(D,Y,B_Y)$ and $D$ cannot be a component of $E$. 
On the other hand, if $D$ is not exceptional$/Y$ but exceptional$/X$, then by definition of log minimal models,  $a(D,X,B)\le a(D,Y,B_Y)=0$ hence  $a(D,X,B)=0$ which again shows that 
$D$ cannot be a component of $E$. Therefore, $E$ is exceptional$/Y$.\\  

\emph{Step 3.} 
Let $B_W$ be the birational transform of $B$ plus the reduced exceptional divisor of $f$, 
and let $C_W$ be the 
birational transform of $C$ on $W$. Pick a sufficiently small $\delta\ge 0$.
Take a general ample$/Z$ 
divisor $L$ so that $K_W+B_W+\delta C_W+L$ is dlt and nef$/Z$.  
Since $(X/Z,B)$ is lc, 
$$
E':=K_W+B_W-f^*(K_X+B)=\sum_D a(D,X,B)D\ge 0
$$
where $D$ runs over the prime exceptional$/X$ divisors on $W$. So,
$$
K_W+B_W+\delta C_W=f^*(K_X+B)+E'+\delta C_W=g^*(K_Y+B_Y)+E+E'+\delta C_W
$$ 
Moreover, $E'$ is also exceptional$/Y$ because for any prime divisor $D$ on $Y$ which is exceptional$/X$, 
$a(D,Y,B_Y)=a(D,X,B)=0$ hence $D$ cannot be a component of $E'$. 

On the other hand, since $Y$ is $\Q$-factorial, there are 
exceptional$/Y$ $\R$-divisors $F, F'$ on $W$ such that  $C_W+F\equiv 0/Y$ and $L+F'\equiv 0/Y$. 
Now run the LMMP$/Y$ on $K_W+B_W+\delta C_W$ with scaling of $L$ 
which is the same as the LMMP$/Y$ on $E+E'+\delta C_W$ with scaling of $L$. Let $\lambda_i'$ and 
$\lambda'=\lim_{i\to \infty} \lambda_i'$ be the 
corresponding numbers. 
If $\lambda'>0$, then by Step 1 the LMMP terminates since $L$ is ample$/Z$. Since $W\to Y$ is birational, 
the LMMP terminates only when $\lambda_i'=0$ for some $i$ which implies that $\lambda'=0$, a contradiction.
Thus, $\lambda'=0$. On some model $V$ in the process of the LMMP, 
the pushdown of $K_W+B_W+\delta C_W+\lambda_i'L$, say 
\begin{equation*}
\begin{split}
K_V+B_V+\delta C_V+\lambda_i'L_V &\\ 
& \equiv E_V+E_V'+\delta C_V+\lambda_i'L_V\\ 
& \equiv  E_V+E_V'-\delta F_V-\lambda_i'F_V'/Y
\end{split}
\end{equation*}
is nef$/Y$. Applying the negativity lemma over $Y$ shows that $E_V+E_V'-\delta F_V-\lambda_i'F_V'\le 0$. But 
if $i\gg 0$, then $E_V+E_V'\le 0$  because 
$\lambda_i'$ and $\delta$ are sufficiently small. Therefore, $E_V=E_V'=0$ as $E$ and $E'$ are effective.\\

\emph{Step 4.} 
We prove that $\phi\colon X\bir Y$ does not contract any divisors. Assume otherwise and let 
$D$ be a prime divisor on $X$ contracted by $\phi$. Then $D^\sim$ the birational transform of 
$D$ on $W$ is a component of $E$ because by definition of log minimal models $a(D,X,B)<a(D,Y,B_Y)$.
Now, in step 3 take $\delta=0$. The LMMP contracts $D^\sim$  
since $D^\sim$ is a component of $E$ and $E$ is contracted. But this is not possible 
because  $K_{X}+B$ is a limit of movable$/Z$ $\R$-divisors and $D^\sim$ is not a component of 
$E'$ so the pushdown of $K_W+B_W=f^*(K_{X}+B)+E'$ cannot negatively intersect a general curve on $D^\sim/Y$. 
Thus $\phi$ does not contract divisors, in particular, any prime divisor on $W$ which is exceptional$/Y$ 
is also exceptional$/X$. 
Though $\phi$ does not contract divisors but $\phi^{-1}$ might contract divisors. 
The prime divisors contracted by $\phi^{-1}$ appear on $W$.\\

\emph{Step 5.} 
Now take $\delta>0$ in step 3 which is sufficiently small by assumptions.  
As mentioned, we arrive at a model $Y':=V$ on which $E_{Y'}+E_{Y'}'=0$ where 
$E_{Y'}$ and $E'_{Y'}$ are the birational 
transforms of $E$ and $E'$ on $Y'$, respectively.
In view of 
$$
K_{Y'}+B_{Y'}\equiv E_{Y'}+E'_{Y'}=0/Y
$$  
we deduce that $(Y'/Z,B_{Y'})$ is a $\Q$-factorial dlt blowup of $(Y/Z,B_Y)$. 
Moreover, by construction $(Y'/Z,B_{Y'}+\delta C_{Y'})$ is also dlt.\\

\emph{Step 6.} As in step 3, 
$$
E'':=K_W+B_W+C_W-f^*(K_X+B+C)=\sum_D a(D,X,B+C)D\ge 0
$$
is exceptional$/X$ where $D$ runs over the prime exceptional$/X$ divisors on $W$.
Since  
$K_W+B_W+C_W\equiv E''/X$, any LMMP$/X$ with scaling of a suitable ample$/X$ divisor $L$
terminates for the same reasons as in step 3. Indeed, 
if $\lambda_i'$ and 
$\lambda'=\lim_{i\to \infty} \lambda_i'$ are the 
corresponding numbers in the LMMP, then  
we may assume $\lambda'=0$ by Step 1; on some model $V$ in the process of the LMMP, 
the pushdown of $K_W+B_W+ C_W+\lambda_i'L$, say 
$K_V+B_V+ C_V+\lambda_i'L_V$ is nef$/X$. But since $W\to X$ is birational, 
there is some $L'\ge 0$ such that $L\sim_\R -L'/X$ hence 
 $E''_V-\lambda_i'L_V'$ is nef$/X$. Now if $i\gg 0$, the negativity lemma implies that 
 $E_V''=0$ hence the LMMP terminates.
 
So, we get a $\Q$-factorial dlt blowup $(X'/Z,B'+C')$ of $(X/Z,B+C)$ 
where $K_{X'}+B'$ is the pullback of $K_{X}+B$ and $C'$ is the pullback of $C$. 
In fact, $X'$ and $X$ are isomorphic outside the lc centres of $(X/Z,B+C)$ because the 
prime exceptional$/X$ divisors on $X'$ are exactly the pushdown of the prime 
exceptional$/X$ divisors $D$ on $W$ with $a(D,X,B+C)=0$, that is, those which are not 
components of $E''$. Since 
$\Supp C$ does not contain any lc centre of $(X/Z,B+C)$ by step 2, 
$(X'/Z,B')$ is a $\Q$-factorial dlt blowup of $(X/Z,B)$ and $C'$ is just
the birational transform of $C$. 
Note that the prime exceptional divisors of  $\phi^{-1}$ 
are not contracted$/X'$ since their log discrepancy with respect to $(X/Z,B)$ are all $0$, and so their birational 
transforms are not components of $E''$.\\

\emph{Step 7.} Remember that $X_1=X$, $B_1=B$, and $C_1=C$. Similarly, put $X_1':=X'$, $B_1':=B'$, and $C_1':=C'$. 
Since $K_{X_1}+B_1+\lambda_1C_1\equiv 0/Z_1$, 
$K_{X_1'}+B_1'+\lambda_1C_1'\equiv 0/Z_1$. Run the LMMP$/Z_1$ on $K_{X_1'}+B_1'$ 
with scaling of some ample$/Z_1$ divisor (which is automatically also an LMMP with scaling of 
$\lambda_1C_1'$). Assume that this LMMP terminates with a log minimal model  $(X_2'/Z_1,B_2')$.
Since $({X_2}/Z_1,B_2)$ is the lc model of $({X_1}/Z_1,B_1)$ and of $({X_1'}/Z_1,B_1')$, 
$X_2'$ maps to $X_2$ and $K_{X_2'}+B_2'$ is the 
pullback of $K_{X_2}+B_2$. Thus, $(X_2'/Z,B_2')$
is a $\Q$-factorial dlt blowup of $(X_2/Z,B_2)$. Since $K_{X_1'}+B_1'+\lambda_1C_1'\equiv 0/Z_1$, 
$K_{X_2'}+B_2'+\lambda_1C_2'\equiv 0/Z_1$ where $C_2'$ is the birational 
transform of $C_1'$ and actually the pullback of $C_2$.
We can continue this process: that is use the fact that 
$K_{X_2}+B_2+\lambda_2C_2\equiv 0/Z_2$ and  
$K_{X_2'}+B_2'+\lambda_2C_2'\equiv 0/Z_2$ to run an LMMP$/Z_2$ on $K_{X_2'}+B_2'$, etc.
Therefore, we can lift the original sequence to a sequence in the 
$\Q$-factorial dlt case assuming that the following statement holds for each $i$:\\

 $(**)$  some LMMP$/Z_i$ on $K_{X_i'}+B_i'$ with scaling of some ample divisor 
 terminates with a log minimal model $(X_{i+1}'/Z_i,B_{i+1}')$.\\

 In this paragraph, we show that we can assume that $(**)$ holds. 
First, note that if $(X/Z,B+C)$ is klt, then $X'\to X$ is 
a small birational morphism and $(X'/Z,B'+C')$ is also klt hence $(**)$ holds 
by [\ref{BCHM}]. Now assume that $(*)$ in Step 1 holds. Then, 
there is a klt $K_{X_i}+\Delta_i\sim_\R K_{X_i}+B_i+\lambda_i C_i/Z$.
 Thus, 
$$
K_{X_i'}+B_i'\sim_\R (K_{X_i'}+B_i')+\epsilon (K_{X_i'}+\Delta_i')=
$$
$$
(1+\epsilon)(K_{X_i'}+\frac{1}{1+\epsilon}B_i'+\frac{\epsilon}{1+\epsilon}\Delta_i')/Z_i
$$
where $K_{X_i'}+\Delta_i'$ is the pullback of $K_{X_i}+\Delta_i$.
If $\epsilon>0$ is sufficiently small, then 
$(X_i'/Z,\frac{1}{1+\epsilon}B_i'+\frac{\epsilon}{1+\epsilon}\Delta_i')$ is klt 
hence $(**)$ again follows from [\ref{BCHM}] under $(*)$. 
Now note that $(**)$ itself is an LMMP under the assumption $(*)$ bearing in mind that 
some $\Q$-factorial dlt blowup of $({X_{i+1}}/Z_i,B_{i+1})$ (which can be constructed as in Step 
6) is a log minimal model of $({X_i'}/Z_i,B_i')$. So, we can assume that 
 $(**)$ holds; otherwise we can replace the original sequence with the one in $(**)$ and repeat 
Steps 2-7 again.

Note that $Y'\bir X'$ does not 
contract divisors: if $D$ is a prime divisor on $Y'$ which is exceptional$/X'$, then it is exceptional$/X$ and 
so it is exceptional$/Y$ by step 6; but then $a(D,Y,B_Y)=0=a(D,X,B)$ and again by step 6 such divisors 
are not contracted$/X'$, a cotradiction. Thus, $(Y'/Z,B_{Y'})$ of step 5 is a log birational model of $(X'/Z,B')$ 
because $B_{Y'}$ is the birational transform of $B'$.
On the other hand, assume that $D$ is a 
prime divisor on $X'$ which is exceptional$/Y'$. Since $X\bir Y$ does not contract divisors by step 4, $D$ is exceptional$/X$. 
In particular, $a(D,X',B')=a(D,X,B)=0$; in this case 
$a(D,Y,B_Y)=a(D,Y',B_{Y'})>0$ 
otherwise $D$ could not be contracted$/Y'$ by the LMMP of step 5 which started on $W$ because 
the birational transform of $D$ would not be a component of $E+E'+\delta C_W$. 
So, $(Y'/Z,B_{Y'})$ is actually a log minimal model of $(X'/Z,B')$.
Therefore, as in step 4, $X'\bir Y'$ does not contract divisors which implies that $X'$ and $Y'$ are 
isomorphic in codimension one. Now replace the old sequence $X_i\bir X_{i+1}/Z_i$
with the new one constructed in the last paragraph and replace $(Y/Z,B_Y)$ with $(Y'/Z,B_{Y'})$.
So, from now on we can assume that $X, X_i$ and $Y$ 
are all isomorphic in codimension one, and that $(X/Z,B+C)$ is $\Q$-factorial dlt.
 In addition, by step 5, we can also assume that 
$(Y/Z,B_{Y}+\delta C_{Y})$ is $\Q$-factorial dlt for some $\delta>0$.\\

\emph{Step 8.}
Let $A\ge 0$ be a reduced divisor on $W$ whose components 
are general ample$/Z$ divisors such that they generate $N^1(W/Z)$. By step 6, $(X_1/Z,B_1+C_1)$ is obtained by running 
a specific LMMP on 
$K_W+B_W+C_W$. Every step of this LMMP is also a step of an LMMP on $K_W+B_W+C_W+\ep A$ 
for any sufficiently small $\ep>0$, in particular, $(X_1/Z,B_1+C_1+\ep A_1)$ is dlt where $A_1$ is the birational transform of $A$. 
For similar reasons, we can choose $\ep$ so that 
$(Y/Z,B_Y+\delta C_Y+\ep A_Y)$ is also dlt. On the other hand, by [\ref{B-II}, Proposition 3.2], 
perhaps after replacing $\delta$ and $\ep$ with 
smaller positive numbers, we may assume that if $0\le \delta'\le \delta$ and $0\le A_Y'\le \epsilon A_Y$, then any 
 LMMP$/Z$ on $K_Y+B_Y+\delta' C_Y+A_Y'$, consists of only a sequence of log flips which are flops with respect to $(Y/Z,B_Y)$. 
Note that $K_Y+B_Y+\delta' C_Y+A_Y'$ is a limit of movable$/Z$ $\R$-divisors for reasons 
similar to those used in Step 2, so no divisor is contracted 
by such an LMMP.\\

\emph{Step 9.} Fix some $i\gg 0$ so that $\lambda_i<\delta$. 
Then, by [\ref{B-II}, Proposition 3.2], there is $0<\tau \ll \ep$ 
such that $(X_i/Z,B_i+\lambda_i C_i+\tau A_i)$ is dlt and such that if we run the LMMP$/Z$ 
on  $K_{X_i}+B_i+\lambda_i C_i+\tau A_i$ with scaling of some ample$/Z$ divisor, then it will be a 
sequence of log flips which would be a sequence of flops with respect to $(X_i/Z,B_i+\lambda_i C_i)$. 
Moreover, since the components of $A_i$ generate $N^1(X_i/Z)$, we can assume that there is an ample$/Z$ $\R$-divisor 
$H\ge 0$ such that $\tau A_i \equiv H+H'/Z$ where $H'\ge 0$ and $(X_i/Z,B_i+\lambda_i C_i+H+H')$
is dlt. Hence the LMMP terminates by (i)
and we get a model $T$ on which both $K_T+B_T+\lambda_i C_T$ 
and $K_T+B_T+\lambda_i C_T+\tau A_T$ are nef$/Z$.  Again since the components of $A_T$ generate $N^1(T/Z)$,  
there is $0\le A'_T\le \tau A_T$ so that $K_T+B_T+\lambda_i C_T+ A_T'$ is ample$/Z$ and $\Supp A_T'=\Supp A_T$.
Now run the LMMP$/Z$ on 
$K_Y+B_Y+\lambda_i C_Y+A'_Y$ with scaling of some ample$/Z$ divisor  where $A_Y'$ 
is the birational tranform of $A_T'$. The LMMP terminates for reasons similar to the above and we end up
with $T$ since $K_T+B_T+\lambda_i C_T+ A_T'$ is ample$/Z$. Moreover, the LMMP consists of only 
log flips which are flops with respect to $(Y/Z,B_Y)$ by Step 8 hence $K_T+B_T$ 
will also be nef$/Z$. So, by replacing $Y$ with $T$ we could assume that $K_Y+B_Y+\lambda_i C_Y$ 
is nef$/Z$. In particular, $K_Y+B_Y+\lambda_j C_Y$ is nef$/Z$ for any $j\ge i$ since $\lambda_j\le \lambda_i$.\\ 

\emph{Step 10.} Pick $j> i$ so that $\lambda_j<\lambda_{j-1}\le \lambda_i$ and let $r\colon U\to X_j$ and $s\colon U\to Y$ 
be a common resolution. Then, we have 

\begin{equation*}
\begin{split}
r^*(K_{X_j}+B_j+\lambda_j C_j) & = s^*(K_Y+B_Y+\lambda_j C_Y)\\ 
r^*(K_{X_j}+B_j) & \gneq s^*(K_Y+B_Y)\\
r^*C_j & \lneq  s^* C_Y
\end{split}
\end{equation*}
where the first equality holds because both $K_{X_j}+B_j+\lambda_j C_j$ and $K_Y+B_Y+\lambda_j C_Y$ 
are nef$/Z$ and $X_j$ and $Y$ are isomorphic in codimension one, the second inequality holds 
because $K_Y+B_Y$ is nef$/Z$ but $K_{X_j}+B_j$ is not nef$/Z$, and the third follows from the 
other two. Now

\begin{equation*}
\begin{split}
r^*(K_{X_j}+B_j+\lambda_{j-1} C_j) & \\
 & = r^*(K_{X_j}+B_j+\lambda_{j} C_j)+ r^*(\lambda_{j-1}-\lambda_j)C_j\\
& \lneq  s^*(K_Y+B_Y+\lambda_{j} C_Y)+s^*(\lambda_{j-1}-\lambda_j)C_Y\\
& =  s^*(K_Y+B_Y+\lambda_{j-1} C_Y)
\end{split}
\end{equation*}

However, since $K_{X_j}+B_j+\lambda_{j-1} C_j$ 
and  $K_Y+B_Y+\lambda_{j-1} C_Y$ are both nef$/Z$, we have 
$$
r^*(K_{X_j}+B_j+\lambda_{j-1} C_j) =s^*(K_Y+B_Y+\lambda_{j-1} C_Y)
$$ 
This is a contradiction and the sequence of log flips terminates as claimed.
$\Box$\\


\vspace{2cm}

\flushleft{DPMMS}, Centre for Mathematical Sciences,\\
Cambridge University,\\
Wilberforce Road,\\
Cambridge, CB3 0WB,\\
UK\\
email: c.birkar@dpmms.cam.ac.uk

\vspace{0.5cm}

Fondation Sciences Math\'ematiques de Paris,\\
IHP, 11 rue Pierre et Marie Curie,\\
75005 Paris,\\
France 


\begin{thebibliography}{99}

\bibitem{}\label{Ambro}  {F. Ambro; {\emph{Quasi-log varieties.}} 
Tr. Mat. Inst. Steklova 240
(2003), Biratsion. Geom. Linein. Sist. Konechno Porozhden-
nye Algebry, 220-239; translation in Proc. Steklov Inst. Math.
2003, no. 1 (240), 214-233. }

\bibitem{}\label{B-II}  {C. Birkar; {\emph{On existence of log minimal models II.}} 
To appear in "J. Reine Angew Math." arXiv:0907.4170v1. }

\bibitem{}\label{B}  {C. Birkar; {\emph{On existence of log minimal models.}}  {Compositio Math.} volume 145 (2009), 1442-1446.}

\bibitem{}\label{BCHM}  {C. Birkar, P. Cascini, C. Hacon, J. M$^c$Kernan; {\emph{Existence of minimal models
for varieties of log general type.}}  J. Amer. Math. Soc. 23 (2010), 405-468. }

\bibitem{}\label{Fujino}  {O. Fujino; {\emph{Fundamental theorems for the log minimal model program.}} arXiv:0909.4445v2.}

\bibitem{}\label{pl-flips} {V.V. Shokurov; {\emph{Prelimiting flips.}}
Proc. Steklov Inst. Math. 240 (2003), 75-213.} 


\bibitem{}\label{log-flips}  {V.V. Shokurov; {\emph{3-fold log flips.}}
With an appendix in English by Yujiro Kawamata.
Russian  Acad. Sci. Izv. Math.  40  (1993),  no. 1, 95--202.}


\end{thebibliography}
\end{document}